\documentclass{amsart}
\usepackage{amsmath,amssymb,amsthm,xcolor}
\theoremstyle{definition}

\theoremstyle{remark}

\numberwithin{equation}{section}

\begin{document}

%\title{A counterexample on sign patterns}

\title{A sign pattern with non-zero elements on the diagonal whose minimal rank realizations are not diagonalizable over the complex numbers}

%    Information for first author
\author{Yaroslav Shitov}

%\address{National Research University Higher School of Economics, 20 Myasnitskaya Ulitsa, Moscow 101000, Russia}
\email{yaroslav-shitov@yandex.ru}

%    \subjclass is required.

\begin{abstract}
The rank of the $9\times 9$ matrix
$$
\left(
\begin{array}{cccc|c|cccc}%{cc|cc|c|cc|cc}
\color{red}{1}&\color{red}{1}&0&0&1&0&0&0&0\\
\color{red}{1}&\color{red}{1}&0&0&0&0&0&0&0\\%\hline
0&0&\color{red}{1}&\color{red}{1}&1&0&0&0&0\\
0&0&\color{red}{1}&\color{red}{1}&0&0&0&0&0\\\hline
0&0&0&0&\color{black}{1}&0&1&0&1\\\hline
0&0&0&0&0&\color{red}{1}&\color{red}{1}&0&0\\
0&0&0&0&0&\color{red}{1}&\color{red}{1}&0&0\\%\hline
0&0&0&0&0&0&0&\color{red}{1}&\color{red}{1}\\
0&0&0&0&0&0&0&\color{red}{1}&\color{red}{1}
\end{array}
\right)
$$
is $6$. If we replace the ones by arbitrary non-zero numbers, we get a matrix $B$ with $\operatorname{rank} B\geqslant6$, and if $\operatorname{rank} B=6$, the $6\times 6$ principal minors of $B$ vanish.
%Whichever non-zero real numbers we put in the place of ones, (i) we get a matrix $B$ of rank at least six, and (ii) if the rank of such $B$ equals six, its principal $6\times 6$ minors vanish.
\end{abstract}

\maketitle

If $\operatorname{rank} B=6$, then the $2\times 2$ red blocks of $B$ are singular. In this case, either the rows or the columns of $B$ corresponding to any of such blocks are collinear, and hence every $2\times 2$ red block can contribute at most one index to a non-singular principal submatrix of $B$. According to Theorem~2.4 in~\cite{manyaut}, the lack of a non-singular  $6\times 6$ minor guarantees that $B$ is not diagonalizable over $\mathbb{C}$. So as the sign pattern of $B$ satisfies the condition as in the title, we get a negative solution to the question that appeared as Problem~2.17 in~\cite{manyaut} and Problem~3.0.15 in~\cite{Zag}.

\end{document}